\documentclass[10pt]{amsart}

\usepackage[utf8]{inputenc}

\usepackage{url}

\usepackage{arydshln}

\usepackage{rotating}

\usepackage[small]{caption}

\usepackage[all]{xy}
\usepackage{pb-diagram, pb-xy}

\usepackage{amssymb}
\usepackage{amsthm}

\usepackage{mathrsfs}

\theoremstyle{plain}
\newtheorem{thm}{Theorem}
\newtheorem{lem}[thm]{Lemma}
\newtheorem{cor}[thm]{Corollary}
\newtheorem{prop}[thm]{Proposition}

\newtheorem{conj}[thm]{Conjecture}

\theoremstyle{definition}
\newtheorem{rem}[thm]{Remark}

\newtheorem*{thm*}{Theorem}

\newcommand{\bbC}{{\mathbb C}}

\newcommand{\bbN}{{\mathbb N}}
\newcommand{\bbP}{{\mathbb P}}
\newcommand{\bbQ}{{\mathbb Q}}

\newcommand{\bbZ}{{\mathbb Z}}

\newcommand{\calA}{{\mathcal A}}

\newcommand{\scrD}{{\mathscr D}}
\newcommand{\scrH}{{\mathscr H}}

\newcommand{\scrT}{{\mathscr T}}

\newcommand{\frakS}{{\mathfrak S}}

\newcommand{\Dbc}{{{D^{\hspace{0.01em}b}_{\hspace{-0.13em}c} \hspace{-0.05em} }}}

\newcommand{\Sp}{\mathit{Sp}}

\newcommand{\Hom}{{\mathit{Hom}}}

\newcommand{\id}{{\mathit{id}}}

\DeclareMathOperator{\Rep}{{Rep}}

\DeclareMathOperator{\CC}{{CC}}

\newcommand{\Atilde}{{\tilde{A}}}
\newcommand{\Dtilde}{{\tilde{D}}}

\begin{document}

\title[]{Cubic threefolds, Fano surfaces and the monodromy of the Gauss map} 
\author{Thomas Kr\"amer}
\address{Centre de Math\'ematiques Laurent Schwartz, \'Ecole Polytechnique, F-91128 Palaiseau cedex (France)}
\email{thomas.kraemer@polytechnique.edu}


\begin{abstract}
The Tannakian formalism allows to attach to any subvariety of an abelian variety an algebraic group in a natural way. The arising groups are closely related to moduli questions such as the Schottky problem, but in general they are still poorly understood. In this note we show that for the theta divisor on the intermediate Jacobian of a cubic threefold, the Tannaka group is exceptional of type $E_6$. This is the first known exceptional case, and it suggests a surprising connection with the monodromy of the Gauss map.
\end{abstract}

\maketitle

\thispagestyle{empty}

\section*{Introduction}

Let $A$ be a complex abelian variety. To any closed subvariety $Z\hookrightarrow A$ one may attach a semisimple complex algebraic group $G_Z$ in a natural way, using the Tannakian formalism of~\cite{KraemerSemiabelian}~\cite{KrWVanishing} as we explain in section~\ref{sec:tannaka}. These Tannaka groups are closely related to classical topics such as the Schottky problem~\cite{KrWSchottky}, intersections of theta divisors~\cite{KraemerE6}, the Torelli theorem~\cite{WeTorelli} and Brill-Noether theory~\cite{WeBN}.~In all previously known examples they are classical groups, and they are simply connected unless $Z$ is a sum of two positive-dimensional subvarieties of $A$. This raises
 several questions: 
(1) Are there examples where $G_Z$ is an exceptional group? 
(2) If~$G_Z$ is not simply connected, then is its simply connected cover always realized as the Tannaka group of some other subvariety? 
(3) What sort of geometric information is encoded in these groups? 
In this note we discuss a family of principally polarized abelian varieties of dimension five which gives a positive answer to the first question, the intermediate Jacobians of smooth cubic threefolds~\cite{CGIntermediateJacobian}:

\begin{thm} \label{thm:E6}
Let $\Theta \subset A$ be the theta divisor on the intermediate Jacobian $A$ of a smooth cubic threefold. Then 
$$G_\Theta \; \cong \; E_6(\bbC)/Z,$$
where $Z$ denotes the center of the simply connected complex algebraic group $E_6(\bbC)$.
\end{thm}

A more precise statement will be given in theorem~\ref{thm:main}. In particular, the universal covering group $E_6(\bbC)$ is realized as the Tannaka group for the Fano surface of lines on the threefold so that the second of the above questions also has a positive answer in this example. Concerning the third question, recall that the Weyl group $W(E_6)$ is the symmetry group of the $27$ lines on a cubic surface. The proof of theorem~\ref{thm:main} will not use this explicitly, indeed it needs surprisingly little geometry and rests on representation theory, see lemma~\ref{lem:rep}. But the result suggests that for any smooth subvariety $Z\hookrightarrow A$ the Weyl group $W(G_Z)$ should contain the monodromy group of a certain Gauss map as a subgroup of small index, see section~\ref{sec:monodromy}; this would explain the relevance of the above Tannaka groups for classical moduli problems.

\medskip

I would like to thank Sebastian Casalaina-Martin, Gavril Farkas, Frank Gounelas and Rainer Weissauer for interesting discussions and helpful comments. I am also very grateful to the organizers of the workshop Classical Algebraic Geometry at Oberwolfach 2014 for the inspiring atmosphere from which this note originated, and to the referee for a careful reading of the manuscript.

\section{The Tannakian framework} \label{sec:tannaka}

Let $D(A)=\Dbc(A, \bbC)$ be the bounded derived category of constructible sheaves of complex vector spaces on a complex abelian variety $A$, and $P(A)\subset D(A)$ its full abelian subcategory of perverse sheaves~\cite{BBD}~\cite{HTT}. The group law $a: A\times A \rightarrow A$ endows these with a rich structure: The convolution product 
\[
 *: \quad P(A) \times P(A) \longrightarrow D(A), \quad \delta_1*\delta_2 = Ra_*(\delta_1\boxtimes \delta_2)
\]
leads to a Tannakian description for the category of perverse sheaves in terms of group representations~\cite{KraemerSemiabelian}~\cite{KrWVanishing}. There is a slight technicality here since in general the convolution of perverse sheaves is no longer perverse. To get around this, let us say a perverse sheaf $\delta\in P(A)$ with hypercohomology $H^\bullet(A, \delta)$ is~{\em negligible} if the Euler characteristic
$$\chi(\delta) \;=\; \sum_{i\in \bbZ} (-1)^i \dim_\bbC H^i(A, \delta)$$ 
vanishes. The negligle perverse sheaves have been classified in~\cite{WeissauerDegeneratePerverse}~\cite[cor.~5.2]{SchnellHolonomic} and may be described explicitly in terms of perverse sheaves on proper abelian quotient varieties of $A$. If $N(A)\subset D(A)$ denotes the full subcategory of all sheaf complexes whose perverse cohomology sheaves are negligible, it has been shown in~\cite{KrWVanishing} that the triangulated quotient category $D(A)/N(A)$ inherits a well-defined convolution product $*$ and that this latter product preserves the full abelian subcategory
$$\overline{P}(A) \;=\; P(A)/(N(A)\cap P(A)) \;\subset\; D(A)/N(A).$$
Furthermore~$\overline{P}(A)$ is a direct limit of Tannakian categories: The subquotients of the convolution powers of any given perverse sheaf $\delta\in \overline{P}(A)$ generate a full abelian subcategory $\langle \delta \rangle \subset \overline{P}(A)$ which is stable under the convolution product and admits an equivalence
\[
 \omega: \quad \langle \delta \rangle \; \stackrel{\sim}{\longrightarrow} \; \Rep(G)
\]
with the Tannakian category of finite-dimensional complex linear representations of a complex algebraic group $G=G(\delta)$. In particular,
 $\omega(\delta_1 * \delta_2) \cong \omega(\delta_1)\otimes \omega(\delta_2)$ for all $\delta_1, \delta_2 \in \langle \delta \rangle$.~It follows from the constructions that the dimension of the arising representations is the Euler characteristic: $\dim_\bbC \omega(\delta) = \chi(\delta)$.

\medskip

In principle this reduces the study of perverse sheaves on $A$ to the representation theory of algebraic groups, though it may be hard to determine the groups $G(\delta)$ explicitly. 
The most interesting case arises when $\delta=\delta_Z$ is taken to be the perverse intersection cohomology sheaf supported on a closed subvariety $Z\hookrightarrow A$. Often~$Z$ is only determined up to a translation by some point of the abelian variety, but the group $G(\delta_Z)$ depends on the chosen translate. So instead of this group we consider the semisimple group
$ G_Z = [G(\delta_Z)^\circ, G(\delta_Z)^\circ ] $
which is the derived group of the connected component.
This derived group remains unaffected if $Z$ is replaced by a translate. Let $\omega_Z = \omega(\delta_Z)|_{G_Z}\in \Rep(G_Z)$ denote its defining representation.

\section{Intermediate Jacobians} \label{sec:intermediate}

Let $V \subset \bbP^4$ be a smooth cubic threefold. Clemens and Griffiths have shown~\cite{CGIntermediateJacobian} that the Fano surface $S$ of lines on $V$ is smooth with $\chi(\delta_S)=27$, embeds in the intermediate Jacobian 
$$ JV \;=\; \Hom(H^{2,1}(V), \bbC)/H_3(V, \bbZ), $$ 
and that the latter is a principally polarized abelian variety of dimension $g=5$ which admits as a theta divisor the image $\Theta = S-S \subset JV$ of the difference morphism $d: S\times S \to JV, (s,t)\mapsto s-t$.  This theta divisor has an isolated singularity at the origin, and there its projective tangent cone is isomorphic to $V$ by~\cite{BeauvilleSingularites}. One then computes $\chi(\delta_\Theta)=78$, see section~\ref{sec:numerical}.

\medskip

The numbers $27$ and $78$ also arise in representation theory as the dimensions of the smallest irreducible representations of the simply connected complex algebraic group $E_6(\bbC)$. Up to duality they are attained precisely for the $27$-dimensional first fundamental representation $\omega_1$ and the $78$-dimensional adjoint representation $Ad$ of this group on its Lie algebra. Note that the adjoint representation factors over the quotient $E_6(\bbC)/Z$ by the center $Z=Z(E_6(\bbC))$.

\begin{thm} \label{thm:main}
Let $A$ be the intermediate Jacobian of a smooth cubic threefold, and consider the corresponding Fano surface $S\subset A$ and the theta divisor $\Theta = S-S \subset A$ as above. Then 
\[ 
 (G_S, \omega_S) \;\cong\; (E_6(\bbC), \omega_1)
 \quad \textnormal{\em and} \quad
 (G_\Theta, \omega_\Theta) \;\cong\; (E_6(\bbC)/Z, Ad). 
\]
\end{thm}

{\em Proof.}
We have already remarked that $\Theta \subset A$ is the image of the difference morphism $d: S\times S \to A$. The latter is generically finite, so by adjunction it follows from the decomposition theorem that $\delta_\Theta$ occurs as a direct summand of the convolution $Rd_*(\delta_S \boxtimes \delta_S)=\delta_S * \delta_{-S} $. On the Tannakian side this gives an embedding
$$ \omega_\Theta \; \hookrightarrow \; \omega_S \otimes \omega_S^\vee 
\quad \textnormal{for the dual} \quad 
\omega_S^\vee \;=\; \Hom(\omega_S, \bbC), $$
so that
$G_\Theta$ becomes identified with the image of the Tannaka group $G_S$ under the tensor product representation $\omega_S \otimes \omega_S^\vee$. The rest of the proof is an application of representation theory which may be found in lemma~\ref{lem:rep} below. \qed

\begin{rem}
A result of Collino~\cite{CollinoFanoI} says that in the moduli space $\calA_5$ of principally polarized abelian fivefolds, the closure of the locus of intermediate Jacobians of smooth cubic threefolds contains the locus of Jacobians of hyperelliptic curves. See also~\cite{CMCubicThreefoldsII}. For a degeneration of a general intermediate Jacobian into a hyperelliptic one, the Fano surface $S$ degenerates to the Brill-Noether subvariety $W_2$. The latter is by definition the image of the symmetric square of the curve in its Jacobian variety, and it has the Tannaka group $G_{W_2} \cong \Sp_8(\bbC)/\pm 1$ by~\cite[th.~6.1]{KrWSmall}~\cite{WeBN}. By semicontinuity~\cite[sect.~4]{KrWSchottky} this group must be a subquotient of the Tannaka group for the Fano surface, which gives a consistency check for our results and could be used for an alternative proof of theorem~\ref{thm:main}:
$$
\xymatrix@M=0.6em@C=0.8em@R=2.5em{
 & \Sp_8(\bbC) \ar@{^{(}->}[rr] \ar@{->>}[d] && E_6(\bbC) \ar@{=}[r] & G_S \\  
 G_{W_2} \ar@{=}[r] & \Sp_8(\bbC)/\pm 1 &&   & 
} 
$$
\end{rem}

\begin{lem} \label{lem:rep}
If $G$ is a connected semisimple complex algebraic group which admits a faithful irreducible representation $V\in \Rep(G)$ of dimension $27$ such that $V\otimes V^\vee$ contains an irreducible direct summand $W\in \Rep(G)$ of dimension $78$, then up to duality
\[
 G \; \cong \; E_6(\bbC), \quad V \; \cong \; \omega_1 \quad \textnormal{\em and} \quad W \; \cong\; Ad.
\]

\end{lem}

{\em Proof.} To any connected semisimple algebraic group $G$ one may associate its universal covering 
$$\tilde{G} \;=\; G_1 \times \cdots \times G_n \; \twoheadrightarrow \; G, $$
where each $G_i$ is simply connected and simple modulo its center. It then follows that 
$$
 V|_{\tilde{G}} \;\cong\; V_1 \boxtimes \cdots \boxtimes V_n
$$
is an exterior product of irreducible representations $V_i\in \Rep(G_i)$ with finite kernel and dimension $d_i=\dim(V_i)>1$. Having $d_1\cdots d_n=\dim(V)=27$ forces that~$n\leq 3$ and $d_i\in \{3,9,27\}$. Now for complex simple Lie algebras of any Dynkin type, the highest weights of all irreducible representations of a given (small) dimension are easily determined via the result of~\cite{AEV}, see table~\ref{tab:rep} where we denote by $\omega_1, \omega_2, \dots$ the fundamental weights. In dimensions $d_i\in \{3,9,27\}$ we are left with the following representations and their duals:

{
\small
\setlength{\dashlinedash}{0.1pt}
\setlength{\dashlinegap}{1.5pt}
\[
 \begin{array}{|r|l:l|l:l:l|l:l:l:l:l:l:l:l|} \hline
 d_i & \multicolumn{2}{l|}{3} & \multicolumn{3}{l|}{9} & \multicolumn{8}{l|}{27} \\ \hdashline
 G_i & A_1 & A_2 & A_1 & A_8 & B_4 & A_1 & A_2 & A_{26} & B_3 & B_{13} & C_4 & E_6 & G_2 \\ \hdashline
 V_i & 2\omega_1 & \omega_1 & 8\omega_1 & \omega_1 & \omega_1 & 26\omega_1 & 2\omega_1 + 2\omega_2 & \omega_1 & 2\omega_1 & \omega_1 & \omega_2 & \omega_1 & 2\omega_1 \\ \hline
 \end{array}
 \medskip
\]
}

\noindent
In our case each $V_i\otimes V_i^\vee$ must by assumption have an irreducible direct summand~$W_i$ such that
$$
 W|_{\tilde{G}} \;\cong\; W_1 \boxtimes \cdots \boxtimes W_n.
$$
Since $\dim(W)=78$, a direct computation shows that the only possibility is $n=1$ and that the universal covering group $\tilde{G}=G_1$ is isomorphic to $E_6(\bbC)$. \qed

\medskip

\section{The magic number $78$} \label{sec:numerical}

It remains to check that for the theta divisor on the intermediate Jacobian of a smooth cubic threefold we have $\chi(\delta_\Theta)=78$. This is a simple computation, but we include it since the outcome is crucial for our proof of theorem~\ref{thm:main}. We begin with a blowup formula.
Let $A$ be a complex abelian variety, and consider an effective divisor $D\subset A$ whose singular locus $D^\mathit{sg} = \{ p_1, \dots, p_n\}$ consists of finitely many isolated points~$p_i$ of multiplicity $m_i$. Blowing up these finitely many points in $A$ we obtain a Cartesian diagram
\[
\xymatrix@R=2em@C=2em@M=0.5em{
 E \ar@{^{(}->}[r] \ar[d] & \tilde{D} \cup E  \ar@{^{(}->}[r] \ar[d] & \tilde{A} \ar[d]^-\pi \\
 D^\mathit{sg} \ar@{^{(}->}[r] & D \ar@{^{(}->}[r]  & A 
}
\]
where $E$ denotes the exceptional divisor of the blowup and where $\Dtilde$ is the strict transform of the divisor $D$. The restriction $\pi: \Dtilde \rightarrow D$ is the blowup of our original divisor in the singular locus and so the scheme-theoretic intersection $E\cap \Dtilde \subset \Dtilde$ is a Cartier divisor. In particular, if the latter is smooth, then $\Dtilde$ is smooth.

\begin{prop}
If $\Dtilde$ is smooth, then

\[
 (-1)^{g-1} \, \chi(\tilde{D}) \;=\; 
 \deg \, [D]^g   \;-\; 
 \sum_{i=1}^n \, m_i \, \biggl[ \frac{(1+t)(1-t)^g}{1-m_i t} \biggr]_{t^{g-1}} \medskip
\]
where $\deg: H^{2g}(A, \bbZ) \rightarrow \bbZ$ denotes the degree map and where the brackets on the right hand side indicate the coefficient of $t^{g-1}$ in the power series.
\end{prop}

{\em Proof.} Let  $E_i = \pi^{-1}(p_i)\subset \Atilde$ denote the exceptional fibres of the blowup, and consider the fundamental cohomology classes $\eta_i = [E_i]$, $\theta = [\Dtilde]$ in $H^2(\Atilde, \bbQ)$. Since the tangent bundle to $A$ is trivial, the formula for the total Chern class of a blowup in~\cite[ex.~15.4.2]{Fulton} says that $c(\Atilde) = (1+\eta)(1-\eta)^g$ where $\eta = \sum_i \eta_i$. So the conormal sequence together with the projection formula for the embedding $\iota: \Dtilde \hookrightarrow \Atilde$ gives the formula
$
 \iota_* ( c(\Dtilde) ) = (1+\eta)(1-\eta)^g \cdot \theta \cdot (1+\theta)^{-1}
$
for the total Chern class of the proper transform. The Gauss-Bonnet theorem says that
$ \chi(\Dtilde)$ 
is the degree $g$ part of $\iota_* ( c(\Dtilde) )$, so our claim follows by writing
\[
 \theta \;=\; \pi^*[D] \;-\; \sum_{i=1}^n m_i \cdot \eta_i
\]
and using that $\eta_i \cdot \eta_k = 0$ for $i\neq k$ and $\eta_i \cdot \pi^*[D] = 0$ whereas $\deg \eta_i^g = (-1)^{g-1}$ for all $i$. This last intersection number is obtained from the conormal sequence for the smooth divisor $E_i\cong \bbP^{g-1}$ with $\deg c_{g-1}(E_i)=\chi(E_i) = g$. \qed

\medskip

\begin{cor} \label{cor:chi}
The theta divisor $\Theta \subset A$ on the intermediate Jacobian $A=JV$ of a smooth cubic threefold $V$ has $$\chi(\delta_\Theta) = 78.$$
\end{cor}

{\em Proof.} Let $\pi: \tilde{\Theta} \to \Theta$ denote the blowup of the theta divisor in the origin. The fibre $\pi^{-1}(0)$ is isomorphic to our threefold $V$ by~\cite[th.~1]{BeauvilleSingularites}, so base change implies that the stalk cohomology of $R\pi_*(\delta_{\tilde{\Theta}})$ in the origin is  
$
 \scrH^i(R\pi_* (\delta_{\tilde{\Theta}}))_0 \cong  H^{i+4}(V, \bbC)
$
for all $i\in \bbZ$. On the other hand 
$$R\pi_*(\delta_{\tilde{\Theta}}) \; \cong \; \delta_\Theta \oplus \varepsilon$$ 
by the decomposition theorem, where $\varepsilon$ is a skyscraper complex which is supported in the origin and is stable under the Lefschetz operator and its inverse. Since the stalk cohomology of $\delta_\Theta$ vanishes in all degrees $i\geq 0$, a comparison with $H^{i+4}(V, \bbC)$ shows $\varepsilon \cong \delta_0[2]\oplus \delta_0 \oplus \delta_0[-2]$ and hence  $\chi(\delta_{\tilde{\Theta}}) = \chi(\delta_\Theta) + 3$.
Therefore our claim follows from the above proposition which in the special case $n=1$, $m_1 = 3$ shows that $\chi(\delta_{\tilde{\Theta}}) = 5! - 39 = 81=78+3$. \qed

\medskip

\begin{rem} \label{rem:index_formula}
(a) The Euler characteristic of a perverse sheaf has a simple meaning in terms of Gauss maps. To explain this, let $A$ be any complex abelian variety, and put $\Omega=H^0(A, \Omega^1_A)$. For any closed subvariety $Z\hookrightarrow A$, the closure of the conormal bundle to the smooth locus of $Z$ is an irreducible subvariety $\Lambda_Z\hookrightarrow \scrT^*_A = A\times \Omega$ of the cotangent bundle to the abelian variety. We define the corresponding Gauss map to be the composite
\[
 \gamma_Z: \quad \Lambda_Z \;\hookrightarrow\; \scrT^*_A \; = \; A\times \Omega \;\twoheadrightarrow\; \Omega,
\]
and we denote the degree of this generically finite map by $\deg(\gamma_Z)\in \bbN\cup \{0\}$. Now to any perverse sheaf $\delta\in P(A)$ we may attach a regular holonomic $\scrD_A$-module by the Riemann-Hilbert correspondence~\cite{HTT}, and its characteristic cycle is a finite formal sum
\[
 \CC(\delta) \;=\; \sum_{Z\hookrightarrow A} m_Z(\delta) \cdot \Lambda_Z
\]
with coefficients $m_Z(\delta) \in \bbN \cup \{0\}$. Franecki and Kapranov~\cite{FKGauss} have shown that in this situation
\[
 \chi(\delta) \;=\; \sum_{Z\hookrightarrow A} m_Z(\delta) \cdot \deg(\gamma_Z).
\]
If $\delta=\delta_W$ is taken to be the perverse intersection cohomology sheaf supported on a closed subvariety $W\hookrightarrow A$, then $m_W(\delta)=1$, but depending on how bad the singularities of the subvariety are, we may also have $m_Z(\delta) > 0$ for some $Z\hookrightarrow W^\mathit{sg}$ inside the singular locus.

\medskip

(b) The theta divisor $\Theta \subset A$ on the intermediate Jacobian $A$ of a smooth cubic threefold has an isolated singularity at the origin, so
$\CC(\delta_\Theta) = \Lambda_\Theta + m\cdot \Lambda_{\{0\}}$ for some $m \geq 0$.
To compute the multiplicity $m$, recall that by~\cite[sect.~13]{CGIntermediateJacobian} we have a commutative diagram of rational maps
\[
\xymatrix@C=1em@M=0.5em{
S\times S \ar[rr]^-\Psi \ar[drr]_-\Phi && \Theta \ar@{=}[r] \ar@{..>}[d] & \bbP \Lambda_\Theta \ar@{..>}[d]^-{\bbP \gamma_\Theta}\\
 && \bbP^4 \ar@{=}[r] & \bbP \Omega
}
\]
where $\Psi$ has generic degree six and where the generic fibre of $\Phi$ can be identified with the set of $27\cdot 16 = 432$ pairs of skew lines on the smooth cubic surface arising as the generic hyperplane section of the threefold. Note that there is a counting error in loc.~cit.~in the paragraph after (13.7); there all $27\cdot 26$ pairs of lines are considered, but the non-skew pairs are mapped by $\Phi$ to a proper closed subset of~$\bbP^4$ and do not contribute to the generic fibre. With the corrected number we get $\deg(\gamma_\Theta)=432/6=72$, so corollary~\ref{cor:chi} and the formula of Franecki and Kapranov imply
\[
 \CC(\delta_\Theta) \;=\; \Lambda_\Theta \;+\; m\cdot \Lambda_{\{0\}}
 \quad 
 \textnormal{with}
 \quad 
 m\;=\; 78 - 72 \;=\; 6.
\]
The nontrivial contribution of the singular locus may be surprising since irreducible theta divisors have only mild singularities~\cite{EL}.
\end{rem}

\section{The monodromy of the Gauss map} \label{sec:monodromy}

Together with the previously known examples, theorem~\ref{thm:main} suggests the following geometric interpretation for the Tannaka groups $G_Z$ from section~\ref{sec:tannaka}.

\begin{conj} \label{conj:monodromy}
For any smooth closed subvariety $Z$ of an abelian variety~$A$, the Weyl group $W_Z=W(G_Z)$ contains as a subgroup of small index the monodromy group $M_Z$ of the Gauss map 
$
 \gamma_Z: \Lambda_Z \longrightarrow \Omega
$.
\end{conj}

 Here the monodromy group is defined as follows. By~\cite{WeissauerDegenerate}~\cite[cor.~5.2]{SchnellHolonomic} the Gauss map $\gamma_Z$ is dominant unless $Z$ is stable under translations by all points of a non-zero abelian subvariety, in which case~$G(\delta_Z)=\{1\}$. Discarding this trivial case, we may assume that the Gauss map restricts over an open dense subset of $\Omega$ to a finite \'etale cover of degree $\deg(\gamma_Z) = \chi(\delta_Z) > 0$. We define its monodromy group~$M_Z$ to be the Galois group of the Galois hull of this cover, which is isomorphic to the image of the monodromy representation on a general fibre of $\gamma_Z$. 
 
 \medskip

\begin{thm}
We have the following monodromy and Weyl groups:
\begin{enumerate}
 \item If $A$ is the Jacobian variety of a smooth projective curve $C\hookrightarrow A$ of genus~$g$, then 
\[
 M_C \;\cong\; W_C \;\cong\;
 \begin{cases} (\pm 1)^{g-1} \rtimes \frakS_{g-1} & \textnormal{\em if $C$ is hyperelliptic}, \\
  \frakS_{2g-2} & \textnormal{\em otherwise}.
\end{cases}
\]
\item If $A$ is the intermediate Jacobian of a smooth cubic threefold with Fano surface $S\hookrightarrow A$, then
\[
 M_S \;\cong\; W_S \;\cong\; W(E_6).
\]
 \item If $A$ is a general principally polarized abelian variety of dimension $g>2$ with theta divisor $\Theta\hookrightarrow A$, then 
\[
\qquad \qquad
 M_\Theta \;\cong\; (\pm 1)^r_0 \rtimes \frakS_r
 \quad \textnormal{\em and} \quad
 W_\Theta \;\cong\;
 \begin{cases} (\pm 1)^r \rtimes \frakS_r & \textnormal{\em if $2\mid g$}, \\
  (\pm 1)^r_0 \rtimes \frakS_r & \textnormal{\em if $2\nmid g$},
 \end{cases}
\]
where $r=g!/2$ and where $(\pm 1)^r_0 = \{ (\epsilon_1, \dots, \epsilon_r) \in (\pm 1)^r \mid \epsilon_1 \cdots \epsilon_r = 1\}$.
\end{enumerate}
\end{thm}

{\em Proof.} For the Weyl groups this follows from theorem~\ref{thm:main} and \cite[th.~6.1]{KrWSmall}~\cite{KrWSchottky}, so it only remains to discuss the monodromy groups. For the Jacobian varieties in (1) we have $\Omega=H^0(C, \Omega^1_C)$ and
\[
 \Lambda_C \;=\; \bigl\{ (p, \omega) \in C \times \Omega \mid \omega(p)=0 \bigr\} \;\hookrightarrow\; \scrT_A^* \;=\; A\times \Omega.
\]
Let $\iota: C \rightarrow \bbP \Omega^*$ be the canonical map and $\overline{C} = \iota(C)$ its image. With the usual identification of points in a projective space and hyperplanes in the dual space, this gives the following factorization for the projectivized Gauss map:
\[
\xymatrix@C=1.5em@M=0.5em{
 \bbP \Lambda_C = \bigl\{ (p, H) \in C \times \bbP \Omega \mid \iota(p)\in H \bigr\} 
 \ar[r]^-\alpha &
 \bigl\{ (\overline{p}, H) \in \overline{C} \times \bbP \Omega \mid \overline{p}\in H \bigr\} 
 \ar[r]^-\beta
 & \bbP \Omega
}
\]
By the uniform position principle~\cite[lemma on p.~111]{ACGH} the monodromy group of $\beta$ is the symmetric group $\frakS_d$ of degree $d=\deg(\overline{C})$. In the non-hyperelliptic case we have $d=2g-2$ and $\alpha$ is an isomorphism, so we are done. It remains to discuss the hyperelliptic case. In that case $d=g-1$ and $\alpha$ is the quotient by the hyperelliptic involution. So the general fibre of the Gauss map consists of $g-1$ pairs of points that are interchanged under the hyperelliptic involution. The monodromy $M_C$ is then a subgroup of the semidirect product $(\pm 1)^{g-1} \rtimes \frakS_{g-1}$ which surjects onto the quotient $\frakS_{g-1}$ via the induced permutation action on a general fibre of $\beta$. It remains to show that for any of the $g-1$ pairs of points in a general fibre of the Gauss map, the group $M_C$ contains a permutation which interchanges the two points of this pair but fixes all the other points in the fibre.
But this follows from the observation that through any branch point $\overline{p}\in \overline{C}$ of the double cover $C\to \overline{C} \cong \bbP^1$ there exists a hyperplane $H\in \bbP \Omega$ which does not meet any other branch point of this cover and over which furthermore the cover $\beta$ is unramified, meaning that this hyperplane is nowhere tangent to the curve $\overline{C}$. Hence part (1) follows. 

\medskip

For the intermediate Jacobian $A$ of a smooth cubic threefold $V$ in (2), recall that the Fano surface $S\subset A$ parametrizes the lines on the threefold. In fact there exists by~\cite[prop.~6]{BeauvilleSingularites} an embedding  $V\hookrightarrow \bbP \Omega^* = \bbP T_0 (A)$ such that for any $p\in S(\bbC)$ the corresponding line $L_p\subset V$ is identified with the projective tangent space to the Fano surface at that point: 
\[
\xymatrix@M=0.5em@C=2em@R=1em{
L_p \ar@{^{(}->}[r] \ar@{=}[d] & V \ar@{^{(}->}[r] & \bbP \Omega^* \ar@{=}[d]\\
\bbP T_p(S) \ar@{^{(}->}[rr] && \bbP T_p(A)
}
\]
Identifying points $H\in \bbP \Omega$ with hyperplanes $H\subset \bbP \Omega^*$ in the dual projective space, we obtain that
\[
 \bbP \Lambda_S \;=\; \bigl\{ (p, H) \in S\times \bbP \Omega \mid L_p \subset H \cap V \bigr\}.
\]
So the fibre of the Gauss map $\gamma_S$ over a general point $H\in \bbP \Omega$ is identified with the~$27$ lines on the smooth cubic surface $H\cap V$, and the monodromy group $M_S$ is the group of permutations of these $27$ lines which is induced by variations of the hyperplane $H\in \bbP \Omega$. By~\cite[VI.20]{SegreCubic} this group is $W(E_6)$, so (2) follows.

\medskip

For part (3) we take a translate of the theta divisor $\Theta \subset A$ which is stable under the involution $\iota = -\id: A\longrightarrow A$. The projectivized Gauss map factors over the quotient by this involution:
\[
\xymatrix@M=0.5em{
 \bbP \Lambda_\Theta \;=\; \Theta
 \ar[r]^-\alpha
 & \Theta/\langle \iota \rangle 
 \ar[r]^-\beta
 & \bbP \Omega.
}
\]
From~\cite[sect.~1]{DebarreDeuxComposantes} we know that over a general point $H\in \bbP \Omega$ in the branch locus of the generically finite map $\beta$, the fibre $\beta^{-1}(H)$ will consist of precisely $r-1$ distinct points where $r=\deg(\beta)=g!/2$. Furthermore $\alpha$ is \'etale over the complement of the finitely many singular points of $\Theta/\langle \iota \rangle$. Restricting $\beta \circ \alpha$ to a general line $\bbP^1 \hookrightarrow \bbP \Omega$ we get by~\cite[th.~1.1]{FultonLazarsfeldConnectivity} a branched cover 
\[
 \xymatrix@M=0.5em{
 X \;=\; \Theta \times_{\bbP \Omega} \bbP^1
 \ar[r]^-{\alpha'}
 & X/\langle \iota \rangle 
 \ar[r]^-{\beta'}
 & \bbP^1
}
\]
of irreducible curves such that the monodromy group of the cover $\alpha' \circ \beta'$ coincides with the one of $\alpha\circ \beta$ (the inclusion of a general line induces an epimorphism on the fundamental groups of the complements of the branch loci).
For $g>2$ we may assume that the double cover $\alpha'$ is \'etale and that $\beta'$ is only simply ramified, and our claim follows from~\cite[th.~1]{BFIrreducibility} after passing to the normalizations of the respective curves. Note that the situation for $g=2$ is different but covered by (1). \qed

\medskip

\begin{table}
\setlength{\dashlinedash}{0.1pt}
\setlength{\dashlinegap}{1.5pt}
\captionsetup{width=0.88\textwidth}
\footnotesize
\begin{sideways}
\begin{array}[t]{|r|c:c:c:c:c:c:c:c:c:c:c:c:c:c:c:c|} \hline
 d & A_2 & A_3 & A_4 & A_5 & A_6 & A_7 & B_3 & B_4 & C_2 & C_3 & C_4 & D_4 & D_5 & E_6 & F_4 & G_2 \\ \hline
 2 &-&-&-&-&-&-&-&-&-&-&-&-&-&-&-&- \\ \hdashline
 3 &-&-&-&-&-&-&-&-&-&-&-&-&-&-&-&- \\ \hdashline
 4 &-&-&-&-&-&-&-&-&-&-&-&-&-&-&-&- \\ \hdashline
 5 &-&-&-&-&-&-&-&-& \omega_2 &-&-&-&-&-&-&- \\ \hdashline
 6 & 2\omega_1 & \omega_2 &-&-&-&-&-&-&-&-&-&-&-&-&-&- \\ \hdashline
 7 &-&-&-&-&-&-&-&-&-&-&-&-&-&-&-& \omega_1 \\ \hdashline
 8 & \omega_1+\omega_2 &-&-&-&-&-& \omega_3 &-&-&-&-& \omega_3, \omega_4 &-&-&-&- \\ \hdashline
 9 &-&-&-&-&-&-&-&-&-&-&-&-&-&-&-&- \\ \hdashline
 10 & 3\omega_1 & 2\omega_1 & \omega_2 &-&-&-&-&-& 2\omega_1 &-&-&-&-&-&-&- \\ \hline
 11 &-&-&-&-&-&-&-&-&-&-&-&-&-&-&-&- \\ \hdashline
 12 &-&-&-&-&-&-&-&-&-&-&-&-&-&-&-&- \\ \hdashline
 13 &-&-&-&-&-&-&-&-&-&-&-&-&-&-&-&- \\ \hdashline
 14 &-&-&-&-&-&-&-&-& 2\omega_2 & \omega_2, \omega_3 &-&-&-&-&-& \omega_2 \\ \hdashline
 15 & 4\omega_1, 2\omega_1 + \omega_2 & \omega_1+\omega_3 & 2\omega_1 & \omega_2 &-&-&-&-&-&-&-&-&-&-&-&- \\ \hdashline
 16 &-&-&-&-&-&-&-& \omega_4 & \omega_1+\omega_2 &-&-&& \omega_4 &-&-&- \\ \hdashline
 17 &-&-&-&-&-&-&-&-&-&-&-&-&-&-&-&- \\ \hdashline
 18 &-&-&-&-&-&-&-&-&-&-&-&-&-&-&-&- \\ \hdashline
 19 &-&-&-&-&-&-&-&-&-&-&-&-&-&-&-&- \\ \hdashline
 20 &-& 3\omega_1, 2\omega_2, &-& \omega_3 &-&-&-&-& 3\omega_1 &-&-&-&-&-&-&- \\ 
    && \omega_1+\omega_2 &&&&&&&&&&&&&& \\ \hline
 21 & 5\omega_1 &-&-& 2\omega_1 & \omega_2&-& \omega_2 &-&-& 2\omega_1 &-&-&-&-&-&- \\ \hdashline
 22 &-&-&-&-&-&-&-&-&-&-&-&-&-&-&-&- \\ \hdashline
 23 &-&-&-&-&-&-&-&-&-&-&-&-&-&-&-&- \\ \hdashline
 24 & 3\omega_1+\omega_2 &- & \omega_1+\omega_4 &-&-&-&-&-&-&-&-&-&-&-&-&- \\ \hdashline
 25 &-&-&-&-&-&-&-&-&-&-&-&-&-&-&-&- \\ \hdashline
 26 &-&-&-&-&-&-&-&-&-&-&-&-&-&-& \omega_4 &- \\ \hdashline
 27 & 2\omega_1+2\omega_2 &-&-&-&-&-& 2\omega_1 &-&-&-& \omega_2 &-&-& \omega_1 & - & 2\omega_1 \\ \hdashline
 28 & 6\omega_1 &-&-&-& 2\omega_1 & \omega_2 &-&-&-&-&-& \omega_2 &-&-&-&- \\ \hdashline
 29 &-&-&-&-&-&-&-&-&-&-&-&-&-&-&-&- \\ \hdashline
 30 &-&-&-&-&-&-&-&-& 3\omega_2 &-&-&-&-&-&-&- \\ \hline
\end{array}
\end{sideways}
\vspace*{1em}

\caption{Up to duality, all irreducible representations of dimension $d\leq 30$ for the complex simple Lie algebras of type $\neq A_1$, with the exception of the defining representations of the classical Lie algebras of type~$A_n, B_n, C_n, D_n$ with $d=n+1, 2n+1, 2n, 2n$. 
We denote the representations by their highest weights, using the fundamental weights $\omega_1,\omega_2, \dots$ labelled as in~\cite{BourbakiLie}.}
\label{tab:rep}
\end{table}


\bibliographystyle{amsplain}
\bibliography{Bibliography}

\end{document}